\documentclass{article}
\usepackage{amsmath,amssymb}
\newtheorem{theorem}{Theorem}
\newtheorem{lemma}{Lemma}

\newtheorem{definition}{Definition}

\begin{document}

\title{ A Weyl Creation Algebra Approach to the Riemann Hypothesis
\author{George T. Diderrich  \footnote{Cardinal Stritch University, 6801 N.Yates Rd, Milwaukee, WI 53217}${}$
\footnote{Carroll University, 100 N.East av.,Waukesha, WI 53186}
 \\ \footnotesize gdiderrich@wi.rr.com, gtdiderrich@stritch.edu, gdiderri@carrollu.edu
\\ \footnotesize http://sites.google.com/site/gtdiderrich/
\\ \footnotesize AMS Primary 11M55, Secondary 81R15 Keywords: Riemann Hypothesis, zeros, quantum field theory }
\date{2/25/13}}
\maketitle
\begin{abstract}
 We sketch a Weyl creation operator approach to the Riemann Hypothesis; i.e., arithmetic on the Weyl algebras with ergodic theory to transport operators.  We prove that finite Hasse-Dirichlet alternating zeta functions or eta functions can be induced from a product of ``creation operators''. The latter idea is the result of considering complex quantum mechanics with complex time-space related concepts. Then we overview these ideas, with variations, which may provide a pathway that may settle RH; e.g., a Euler Factor analysis to study the quantum behavior of the integers as one pushes up the critical line near a macro zeta function zero.
\end{abstract}
\tableofcontents
\section{Slowing Down the Integral Expression for $\zeta(s)$}
We assume some familiarity with the classical version of Riemann zeta function theory found in Riemann's original paper \cite{R} 
including elementary complex variable theory and the book by H.M. Edwards \cite{E} and \cite{B2}, \cite{Der}, \cite{Pat}. Conrey gives an overview of some of the previous results in the field \cite{Con2} along with the more recent overview by Sarnak \cite{Sar}. 

 First recall the integral expression of the Riemann zeta function
 
	\[ \int^{+\infty}_{0} {\frac{x^{s-1}}{e ^{x} -1}}dx = \Pi(s-1)\zeta(s) \qquad ( s > 1).
\]
	
\textsl{Note: We use the notation that $ s > 1 $ where $ s = \sigma +it $ means that the real part of s is greater than one i.e.  $\Re s >1 $.}

Instead of the usual "`kernel"' in the definition of the Riemann zeta function above of 
	\[e ^{z} -1 
\]we consider first a ``harmonic'' polynomial of conjugates
	\[ H^{*}_{n}(z) = (z^{2}+1)(z^{2}+2^{2})...(z^{2}+n^{2})
\]
from considerations of quantum mechanics \cite{D1}, \cite{Bo}, \cite{Wyl1}, \cite{Fey1}, \cite{Fey2}, quantum oscillator theory and quantum field theory \cite{D2}, \cite{K}, \cite{Dy}, \cite{Mc}, \cite{Z}, \cite{Ez}. This idea is the outgrowth of previous unpublished research on complex variable quantum oscillators \cite{Did2} (c.f. Section 4) which led to the following paradigm.

\subsection{Introduction of the Paradigm: Linear Forms to Creation Operators}
\textbf{Paradigm} \textit{ One can convert a polynomial in $ z $ over the complex numbers $ \mathbf{C}$  into annihilation and creation operators as follows :}

	\[ z + ik \stackrel{}{\rightarrow} \frac{p + ikz}{\sqrt{2}} = a_{(k)}  
\] 
\textit {where p corresponds to the momentum operator; k corresponds to an angular frequency; and z corresponds to the position operator using essentially Dirac's definitions with all variables complex variables and with complex parameters.}

  Actually, Dirac uses the older notation $ \eta_{k}$ for particle creation operators instead of the $ a_k $ or $ a^{*} $ notation.
It suffices to consider creation operators because the destruction operators are just the anti-particle reversed time induction into the macro space-time thread from it's neighboring space-time thread.
	
 Thus the harmonic polynomial above corresponds to a composite family of quantum oscillators with frequencies equal to the whole numbers. This will give all of the pure imaginary poles of the Riemann zeta kernel up to a factor of $ 2\pi $ and up to $n$ while slowing down the growth properties of
  $ e ^{z} -1 $.
 Furthermore, one can now perform arithmetic c.f. Section 7 using Weyl algebras because the index $k$ in $ a_{(k)} $ possesses arithmetic and analytic properties. Thus essentially we have found a finite quantum field theory to analyze zeta functions. 
	
 Hence, the present line of research is to use standard methods along with our new ideas and see what we can learn to track the alignment of the non-trivial critical zeros of the zeta function. 
\subsection{Some Approximate Dirichlet-Hasse Eta Functions Related to RH} 
 We consider integrals of the form
	\[L = L_n(s) = \int^{+\infty}_{0} {\frac{x^{s-1}}{H^{*}_{n}(x)}}dx \qquad ( 0<s<2n)
\] 
and prove that $ H^{*}_{n}(z)$ induces the formation of a Hasse type of eta function with the properties:
	\[
\]
\begin{theorem} The trivial zeros of 

	\[ \zeta^{H^{*}}_{n}(s) = \sum^{n}_{k=1} (-1)^{(k-1)} \binom{2n}{n+k}k^{-s} 
\]
are given by the negative even integers

	\[\zeta^{H^{*}}_{n}(-2m)= 0  \qquad (1 \leq m \leq n-1)
\]
where  
\[ H^{*}_{n}(z) = (z+i)(z-i)(z+2i)(z-2i)...(z+ni)(z-ni).
\] 
F furthermore, $ \zeta^{H^{*}}_{n}(s) $ is an entire function and $ L_n(s) $ is given by

	\[L_n(s) = \frac{\pi}{sin( s\frac{\pi}{2}) } \zeta^{H^{*}}_{n}(-s) \frac{1}{(2n)!}   \qquad ( 0<s<2n)
\]which extends meromorphically throughout the s-plane with possible poles at the even integers except for the
pole-free points at 
	\[  s=2m \qquad  (  1 \le m \le n-1).
\]
The special case of this theorem with $s=1$ and $n=1$ is the well known result 

	\[ \int^{+\infty}_{0} {\frac{1}{x^{2}+1}}dx = \frac{\pi}{2}.
\]
\end{theorem}
We call $ \zeta^{H^{*}}_{n}(s) $ an approximate alternating zeta function or an approximate eta function corresponding to the polynomial $H^{*}$. This is of interest because the Dirichlet eta function or alternating zeta function

	\[\eta(s)= \sum^{\infty}_{k=1}\frac{(-1)^{k-1}}{k^{s}} = (1-2^{1-s})\zeta(s) \qquad ( s > 0 )
\] 
can be expressed in a slightly different notation using a result of Hasse \cite{Wolf} as follows 

\begin{gather*}
\eta(s) =  \sum^{\infty}_{n=0} \frac{1}{2^{n+1}}\zeta^{H}_{n}(s) \qquad  ( s > 0 )  \\
\eta(s) =   \sum^{\infty}_{n=0} \frac{1}{2^{n+1}}\eta_{n}(s)  \qquad where \\                                                           
\zeta^{H}_{n}(s) = \eta_{n}(s)= \sum^{n}_{k=0}(-1)^{k} \binom{n}{k} (k+1)^{-s}  \qquad and \\
\ H_{n}(z) = (z+1)(z+2)...(z+n+1).
\end{gather*}

In fact we now find a direct connection to Hasse's eta function induced from $ H_{n}(z) $ with the properties:
 
\begin{theorem}
  The trivial zeros of
	\[ \zeta^{H}_{n}(s) = \eta_{n}(s) =\sum^{n}_{k=0} (-1)^{(k)} \binom{n}{k}(k+1)^{-s} 
\] 
are given by negative integers and zero in the range

	\[\zeta^{H}_{n}(-m)=\eta_{n}(-m)= 0  \qquad (0 \leq m \leq n-1)
\] 
where 
\[ H_{n}(z) = (z+1)(z+2)...(z+n+1).
\] 
Furthermore, $ \zeta^{H}_{n}(s) $ is an entire function and $ L_n(s) $ is given by
	\[L = L_n(s) = \int^{+\infty}_{0} {\frac{x^{s-1}}{H_{n}(x)}}dx \qquad ( 0 < s < n+1)   
\]
with
	\[L_n(s) = \frac{\pi}{sin( s \pi )} \zeta^{H}_{n}(1-s) \frac{1}{(n)!}   \qquad ( 0 < s < n+1)
\]which extends meromorphically throughout the s-plane with possible poles at the integers except for the
pole-free points at 

	\[  s=m \qquad  (  0 \le m \le n-1).
\] 
Note that the special case of this theorem with $ n=0 $ is the well known result

	\[ L_{0}(s) = \int^{+\infty}_{0} {\frac{x^{s-1}}{x+1}}dx  = \frac{\pi}{sin( s \pi )} \qquad (0 < s < 1)
.\] 
\end{theorem}

 This shows that there is some virtue in exploring the local approximate $ \eta $ functions to reveal clues about the global behavior of the zeta function in the critical strip. In fact, we speculate that for $ n=0,1,2..$ that the critical line of each of the approximate eta functions $  \eta_{n }(s)$ is on the line $ \frac{n+1}{2} +it $, in a weaker topology, because the region of convergence of $L_{n}(s)$  is $0 < s < n+1 $.
 Moreover, using Sections 5 and 6 as a blueprint, we outline an approach to RH based on the above results and ideas. We plan on studying what we call "`proto-zeros"' to give information about the non-trivial zeros of the zeta function.

 A proto-zero for $\eta_{k}(s)$ is defined to be for $s=\sigma + it $ with fixed $ \sigma$ the smallest in absolute value in a neighborhood of $ it $ of the function. For example, the proto-zeros of 
$\eta_{1}(s)= 1-\frac{1}{2^s} $ for $ s= \sigma + it $ with $ 0 < s < 1$ i.e. occur whenever $ t=\frac{n}{log2/2\pi} $ with $n$ an integer.
 We will also be looking for local functional equations. 
 
\subsection{Riemann's Functional Equation for the Zeta Function.}
 It is interesting, along this line of thinking 
to record the functional equation of Riemann in the following form \cite{E}, \cite{Ca}

	\[ \frac{\zeta(s)}{\zeta(1-s)}=(2 \pi)^{s-1} 2 sin(\frac{\pi s}{2}) \Pi(-s) 
\]

from which it follows fairly readily that the zeta function $\zeta(s) $ has trivial zeros at the negative even integers. Also the above left hand side suggests the symmetry
	\[s \stackrel{}{\rightarrow}  1-s
\] which maps  the critical region
	\[ 0 < s < 1
\] into itself and keeps invariant the critical line
	\[s= \frac{1}{2} +it
\] 
where presumably all of the non-trivial zeros lie. The latter is a statement of the Riemann Hypothesis.
Note that the right hand side of the functional equation contains product factors of the type

	\[ 1+ \frac{s}{\alpha} .
\]

This suggests that more might be done with this observation in attacking the RH c.f. Section 7 on Euler factors approach. Also our results in Theorem 1 and Theorem 2 seem to echo some sort of functional equation connections to Riemann's functional equation but we have not been able to find a direct relationship so far.

\section{Proof of Theorem 1}
We follow almost in parallel the same techniques that may be found in H.M. Edwards \cite{E} for elementary Riemann zeta function theory with the major difference that an explicit functional equation is apparently not manifest but implied from our analysis. 
In our case, we develop the proof for our theorem as follows. First we choose the branch cut corresponding to the non-negative
reals for $log(z)$, which is our notation for $ln(z)$, instead of the usual non-positive reals. Hence the phases for Sheet(0) correspond to the phase band 
  \[ B_{0} = [0,2\pi)
\]
in the counter clock wise positive sense.
Consider the "`branch cut"' contour integral

	\[L_{\gamma}= L_{\gamma n}(s) =  \oint_{\gamma(0)} {\frac{z^{s-1}}{H^{*}_{n}(z)}}dz  \qquad ( s < 2n).
\]
This clearly defines a holomorphic or analytic function using standard methods i.e., uniform convergence on compact domains because $ 2n - (s-1) > 1 $. It is understood that the contour $\gamma(0) $ is the path that starts at  $+\infty $ 
above the non-negative axis and then circles around  zero $ 0 $ in a positive sense back underneath the non-negative axis excluding all other discontinuities of the argument of the integral back to  $+\infty $ .  Define the following path $\gamma $
(c.f. Figure 1)
    \[ \gamma = \gamma(0) + \gamma(\infty) + \gamma(R).
\]
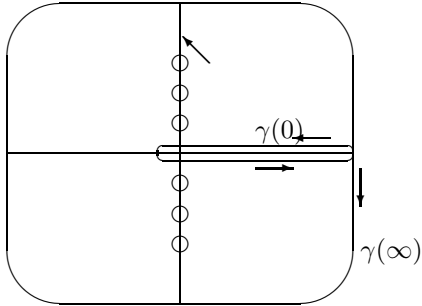
\begin{figure}[h]
\centering 
\setlength{\unitlength}{2cm}
\begin{picture}(3,3)
\put(1.5,1.5){\line(0,1){1.}}
\put(1.5,1.5){\line(0,-1){1}}

\put(1.5,1.5){\line(1,0){1.15}}
\put(1.5,1.5){\line(-1,0){1.15}}

\put(1.5,1.7){\circle{.1}}
\put(1.5,1.3){\circle{.1}}

\put(1.5,1.9){\circle{.1}}
\put(1.5,1.1){\circle{.1}}

\put(1.5,2.1){\circle{.1}}
\put(1.5,.9){\circle{.1}}
\put(1.7,2.1){\vector(-1,1){.18}}

\put(2.,1.5){\oval(1.3,.1)} 
\put(2.5,1.6){\vector(-1,0){.25}}
\put(2.,1.4){\vector(1,0){.25}}
\put(2.,1.6){$ \gamma(0) $}

\put(1.5,1.5){\oval(2.3,2)}
\put(2.7,1.4){\vector(0,-1){.25}}
\put(2.7,.8){$ \gamma( \infty)  $}
\end{picture}
\caption{The path $\gamma $.} 
\label{fig:path} 
\end{figure}
The path $\gamma(\infty)$ is the the path moving downward in a negative sense from  $+\infty $ around the entire z-plane
at  $\infty $ back to the beginning just above the  $+\infty $ point where the path $\gamma(0)$ starts and ends. The path
$\gamma(R)$ refers to all of the little positive sense circles around each of the poles

\[ z_k = ik \qquad ( k=1,2...n,-1,-2,..-n)
\]
where it is understood that the entry and exit lines to each little circle from $ \infty $ and back out to $ \infty $ cancels.
Hence the path $\gamma $ encompasses a single valued holomorphic function defined on Sheet(0) for $ log(z) $ with phase
band $ B_{0} = [0,2\pi) $. So the path integral of this must equal 0 by Cauchy's Theorem:

	\[\oint_{\gamma} = 0 = \oint_{\gamma(0)}+ \oint_{\gamma(\infty)}+ \oint_{\gamma(R)}
\]
Since $ s < 2n $ we have
	\[ \oint_{\gamma(\infty)} = 0
\] 
hence
	\[\oint_{\gamma(0)} = -\oint_{\gamma(R)}
\]
We use the following definitions
	\[ z^{s-1} =\exp (  (s-1) ( log|z| + i \phi(z)       )                      )
\] 
and the evaluation of $ i^{s-1} $ is given by

	\[	E = e^{  \frac{(s-1)\pi}{2}      }.
\]
Using our previous definitions we can evaluate
	\[\oint_{\gamma(0)} = L_{\gamma} = L_{\gamma n}(s) = -L + E^{4}L \qquad( 0<s<2n)
\] 
since the "`phase charge"' on the outbound lower path is $2\pi$.
Also note that

	\[\oint_{\gamma(R)} = 2\pi i \sum r_{k}=2\pi i R \qquad (1\le |k|\le n)
\]
gives the sum of the residues at the poles $ z_{k} = ik $.
Our next task is to compute the residues $r_{k}$ for $k > 0$. Recall that

	\[\lim_{z\to ik} (z-ik) \frac{z^{s-1}}{H_{n}(z)} = r_{k}
\] 
where
	\[ H_{n}(z) = (z+i)(z-i)(z+2i)(z-2i)...(z+in)(z-in) .
\]
To organize the calculations, write

\begin{gather*}
A(z) = (z-i1)(z-2i)...(z-(k-1)i)\\
B(z) = (z-i(k+1))(z-(k+2)i)...(z-ni) \\
C(z) = (z+i)(z+2i)...(z+ni).
\end{gather*}
Therefore,
\begin{gather*}
A(ik) = ((k-1)i)((k-2)i)...(i)= (k-1)! (i)^{k-1}\\
B(ik) = (-i)(-2i)...(-(n-k)i) = (n-k)!(-i)^{n-k} \\
C(ik) = ((k+1)i)((k+2)i)...((n+k)i) =\frac{(n+k)!}{k!}(i)^{n}.
\end{gather*}
Hence,
	\[r_{k} = (ik)^{s-1}\frac{k!}{((k-1)!)(n-k)!(n+k)!(i)^{k-1}(-i)^{n-k}(i)^{n}}
.\]
Simplifying where recall $ 	E = e^{ (s-1)(\frac{1 \pi}{2} ) } $,
	\[r_{k} = E\frac{k^{s}}{(n-k)!(n+k)!(i)^{2n-1}(-1)^{n-k}}
.\]
Finally,
	\[(i)^{2n-1}(-1)^{n-k} = (-1)^{n}(i)^{-1} (-1)^{n-k} = (-1)^{2n-k}(-i) =-i(-1)^{k} .
\]
thus we obtain,
	\[r_{k} = +i E\frac{k^{s} (-1)^{k}}{(n-k)!(n+k)!} \qquad (k>0).
\]
A similar calculation for $k<0$ gives a "`conjugate"' 
	\[r_{k} = -i E^{3}\frac{k^{s} (-1)^{k}}{(n-k)!(n+k)!} \qquad (k<0)
\]
because on Sheet(0) with phase band $B =[0,2\pi)$ we have

	\[ (-i)^{s-1} = e^{ (s-1)(\frac{3 \pi}{2} ) } = E^{3}.
\]
Thus the sum of the residues are
\begin{gather*}
R =  \sum {r_{k}} = i(E-E^{3}) \sum_{k=1}^{n} \frac{k^{s}(-1)^{k} }{ (n+k)!(n-k)! } \\
R = -i(E-E^{3}) \frac{1}{(2n)!} \zeta^{H^{*}}_{n}(-s)    \qquad where          \\
\zeta^{H^{*}}_{n}(s) = \sum^{n}_{k=1} { (-1)^{k-1}} \binom{2n}{n+k}k^{-s}.
\end{gather*}
Furthermore note that $ E-E^{3} = E(1-E^{2}) $ and that 
	\[ 1-E^{2} = 0
\] 
if and only if for the isolated points $ s = 2m +1 $ where m is an integer because  $ \frac{\pi}{2}(s-1) = \pi mod(\pi) $.

Using the previously developed calculations starting from the result that the branch cut contour integral equals
the negative of the residues integral 

	\[\oint_{\gamma(0)} = -\oint_{\gamma(R)}
\] we find

	\[-(1-E^{4}) L  = -2\pi i R 
\]
	\[(1+E^{2})(1-E^{2}) L =  2\pi (i)(-i) E(1-E^{2}) \zeta^{H^{*}}_{n}(-s) \frac{1}{(2n)!} .
\]
 We can cancel $1-E^{2} $ on both sides because this factor has isolated zeros 
only at the odd integers $  1 \le s= 1+2m \le 2n-1 $ therefore the other factor must be zero throughout the specified domain, so we obtain

	\[ (1+E^{2}) L = 2\pi  E \zeta^{H^{*}}_{n}(-s) \frac{1}{(2n)!}
\] thus
	\[ cos((s-1)\frac{\pi}{2}) L = sin( s\frac{\pi}{2}) L = \pi \zeta^{H^{*}}_{n}(-s) \frac{1}{(2n)!} .
\] We also obtained the same result above by working with a path $\gamma $ just above the real axis.
Now by our working assumptions of $ 0 < s < 2n $ we have that

  \[ L = L_n(s) = \int^{+\infty}_{0} {\frac{x^{s-1}}{H^{*}_n(x)}}dx 
   \] 
 is convergent, and that $ \zeta^{H^{*}}_{n}(s) $ is an entire function because it is a finite sum of exponentials. Now assuming even integers $ s = 2m $ 
, we arrive, using the above equation, through essentially zero-pole cancellation,
 
	\[  \zeta^{H^{*}}_{n}(-2m) = 0 \qquad ( 1\le m \le n-1).
\]
Also we have,

	\[L_n(s) = \int^{+\infty}_{0} {\frac{x^{s-1}}{H^{*}_n(x)}}dx = \frac{\pi}{sin( s\frac{\pi}{2}) } \zeta^{H^{*}}_{n}(-s) \frac{1}{(2n)!} \qquad ( 0 < s < 2n)
.\]
which completes the proof of the theorem.

\section{ Proof of Theorem 2}
The proof of Theorem 2 is very similar to the previous Theorem 1; hence, we present only the highlights. 

To begin, the integrand in the corresponding integral is now given by,

	\[\frac{z^{s-1}}{H_{n}(z)}
\]
where
	\[H_{n}(z) = (z+1)(z+2)...(z+n+1) 
\]
is now a product of creation operators $ a_{k}  $  for  $ k = 1,...n+1 $ using our QFT paradigm.
Let 
	\[E = e^{ \frac{i\pi}{2} (s-1)   }
\] 
then, proceeding as before, the residues are given by

	\[ r_{-k} = \frac{E^{2} k^{s-1}   }{(-(k-1))(-(k-2))..(-1)(1)(2)..(n+1-k)} \qquad (1\le k \le n+1)
\]
then changing notation where

	\[ k' = k-1
\]
then gives

\[r_{ -(k'+1)} =  \frac{E^{2} (k'+1)^{s-1} (-1)^{k'}  }{ (k')! (n-k')!  } \qquad( k'= 0,1,..,n).
\]
Therefore the sum of the residues is

	\[ R = \frac{E^{2}}{n!} \sum^{n}_{k'=0} (-1)^{k'} \binom{n}{k'} (k'+1)^{s-1}
\]
hence

	\[R = \frac{E^{2}}{n!} \eta_{n}(1-s)
\]
where 
	\[ \eta_{n}(s) = \zeta^{H}_{n}(s) = \sum^{n}_{k=0} (-1)^{k} \binom{n}{k} (k+1)^{-s}
\]
is Hasse's local eta function in the alternating zeta function mentioned in the introduction.

 Continuing with the calculations, using the fact that the sum of the little residue integrals plus the branch cut integral must
add to zero, gives the relationship,

	\[ -L + E^{4}L  = -2\pi i R
\]
where
	\[ L = \int^{+\infty}_{0} {\frac{x^{s-1}}{H_{n}(x)}}dx \qquad ( 0 < s < n+1).
\]
Thus we find,

	\[ L = 2\pi i \frac{E^{2}}{1-E^{4}}\frac{\eta_{n}(1-s)}{n!}
\] 
which upon simplification gives

\begin{gather*}
L = \frac{\pi}{-sin( (s-1)\pi   )}\frac{\eta_{n}(1-s)}{n!} \qquad or \\
L = \frac{\pi}{sin(s\pi)}\frac{\eta_{n}(1-s)}{n!} .
\end{gather*}
From the above, it is clear that since
\begin{gather*}
\sin(\pi s) = 0  \mbox{ the the poles occur at the integers,} \\
s = 0,\pm 1,\pm 2,.. \mbox{  but the interval of convergence L of is    }\\
0 < s <  n+1 \mbox{, thus} \\
0 > -s > -n-1 \\
1 > 1-s > -n \\
-m = 1-s \\
1 > - m > -n \\
0 \le m \le n-1 \\
 \mbox{ therefore using zero pole cancellation forces}\\
\eta_{n} (-m) = 0 \mbox{   giving the zeros in the above range.  }
\end{gather*}
The remaining details are easily filled in and we have proved the theorem.

 Next, we make some preparatory observations on complex quantum oscillator theory which gives the motivation behind our investigation into RH using QFT ideas. A portion of these results are from my unpublished research \cite{Did2} with inspiration from Lubkin \cite{L1}. 

\section {Quantum Complex Oscillators: The Origin of the Paradigm}
This section requires a little bit of background in elementary quantum mechanics and the standard quantum harmonic oscillator theory \cite{D1}, \cite{Bo}. We will extend the ordinary quantum oscillator Hamiltonian $ H = \frac{1}{2m}(p^{2} + m^{2}\omega^{2}x^{2}) $  to a complex oscillator with complex frequency $\omega $; with complex mass $m$; with complex $ \hbar $; and complex spatial position $ z $; complex momentum $p$; and complex time $ s $.

 We prove that a complex oscillator can be quantized. This will be demonstrated by a symbolic and algebraic method and also by an analytic method.
One should note that Becker \cite{Be1} and Youssef \cite{Yo},\cite{Yo1} others have studied some of these ideas before the author became aware of their publications.
	
 Also, we have since learned that D-module theory \cite{Hartillo},\cite{Ginz} essentially captures the algebraic aspects of these constructions with the introduction of Weyl algebra machinery over the complex numbers $\mathbf{C}=k$. 

 A Weyl algebra is an associative unital algebra given by the free algebra $k[b_{(1)},...b_{(n)},a_{(1)},...a_{(n)}]$ modulo two sided commutator conditions of the form $[b,a] - 1 $ . Specifically,
 
	\[ A(k)_{n} = k[b_{(1)},...,b_{(n)},a_{(1)},...,a_{(n)}]
\]
satisfying commutator conditions with all operators commuting except for

	\[   [b_{(k)},a_{(k)}] = 1 = b_{(k)}a_{(k)} - a_{(k)},b_{(k)}
\] 
where
	\[ b_{(k)} = a_{(-k)}
\]
for all integers $ k = 1,2,..,n$.
 
  The unit operator $1$ commutes with all the $a's$ and $b's$. Also later in our constructions, we allow the indices $k$ to be any complex number. Also note that there is a scaling of the a's and b's operators to obtain the $1$ operator.
 
 The difference in our algebraic approach is to emphasize the introduction of creation operators through the paradigm (c.f. Section 1)

	\[z+ik \stackrel{}{\rightarrow}  \frac{p +ikz}{\sqrt{2}}
\]
and to create "particles" by hitting the unit operator with powers of a fixed creation operator. This has a similar flavor to Nakajima's creation operator methods in Hilbert scheme theory \cite{Nak}.  The anti-particles correspond to negative integer particles. Also we can emulate usual quantum mechanics with kets and bras inside the algebra to compute matrix elements. Hence we view a Weyl algebra as a finite QFT quantum field theory.

 The following provides the original rendition of our thinking which has served the foundation for our current research activities which we believe has intrinsic value. 

\subsection{ Preliminary Quantum Oscillator Calculations}
We start by considering 

	\[ H = \frac{1}{2}(p^{2} + \omega^{2}z^{2}).
\] 
which amounts to assuming $m=1$ and $\hbar =1 $. This can be accomplished by using a series of scaling transformations.

 We have to be careful in specifying the phase of $ \omega $ because we need a rule to distinguish
between creation and destruction operators and hence the given arrow of time. We will see this later.

 Notice that H above is symmetric with the interchange
of $\omega$  and   $-\omega$. For the moment, we ignore this ambiguity and proceed with the following definitions.

 Introducing Schrodinger's prescription,

	\[    p = -i\hbar \partial_{z} = -i\partial_{z}.
\] 
leads to the bosonic commutator quantum condition
	\[ [p,z] = pz-zp = -i
\]
because operating the commutator on the left with an arbitrary spatial wave $\psi(z) $ gives
\begin{gather*}
 (pz -zp)\psi = p(z\psi) - z(p\psi) \\
 [p,z]\psi = -i\psi + z(p\psi) -z(p\psi) = -i\psi.
\end{gather*}
Now classically, at the macro level, factorize $H$ above

	\[ H = \frac{1}{\sqrt{2}}(p + i\omega z)\frac{1}{\sqrt{2}}(p -i\omega z)
\] 
which suggests the following definitions

\begin{gather*}
a = a_{(\omega)} = \frac{1}{\sqrt{2|\omega| }}(p + i\omega z) \\
b = a_{(-\omega)} = \frac{1}{\sqrt{2|\omega|}}(p - i\omega z)  \qquad with \\
\omega =  |\omega  | u_{\omega} = |\omega  | u  \qquad where \\
 u = e^{i \phi}  \mbox {  is a phase factor for the frequency.}
\end{gather*} 
with possible phase ambiguity in the factor $ \sqrt{\omega} $.

This gives the following commutator

	\[  [b,a] = ba - ab = u
\] 
because
\begin{gather*}
ab = \frac{H}{|\omega|} + \frac{i\omega}{2|\omega|}(zp-pz  ) \\
ab = \frac{H}{|\omega|} + u\frac{i i}{2} \\
ab = \frac{H}{|\omega|} - u\frac{1}{2} .
\end{gather*}
	\[
\]
Note: The usual notation for the commutator is $ [a,a^{*}] = 1 $ however, we want to avoid using $C^{*}$ concepts in our deliberations.

 Hence we see that $ ab $ drops the intrinsic energy by one half of a quantum under the tacit assumption
of the given orientation of "`time"'. Similarly, we obtain
	\[ ba = \frac{H}{|\omega|} + u\frac{1}{2}
\]
and that this gives an increase of intrinsic energy by one half quantum.

 Subtracting the two "`number"' operators gives the quantum condition for particles

	\[  [b,a] = ba - ab = u
\]
and adding the two equations together gives H as a symmetrical sum of number operators,
 
	\[  H = \frac{|\omega|}{2} (ba + ab ).
\]
Sometimes we write this in the following form to show the decomposition of the oscillator into a magnitude part times a phase part:

	\[H_{(\omega)} = |\omega| H_{(u)} = |\omega| u H_{(1)} =\omega H_{(1)}
\]
where the subscript defines the type of frequency in an obvious notation.
 Next we bring in "`time". The Hamiltonian equations of motion for a dynamical variable $v$ using Poisson Brackets are

	\[\dot{v} = [v,H]_{PB} .
\]
Following Dirac \cite{D2} this becomes, using commutators in a quantum mechanical description, Heisenberg's equation of motion

	\[   i\hbar \dot{v} = [v,H] = vH - Hv.
\]
Next we examine the law of motion for $a$ assuming a complex frequency $ \omega $ to see the impact on "`time"' and the emergence
of anti-particle concepts. Also we carry out enough of the commutator calculations for the reader to see what the impact of the phase factor ``u'' is on the results.
 
Write assuming again $ \hbar = 1 $ then

	\[     i \dot{a} = [a,H] = \frac{a|\omega|}{2} (ab + ba) - \frac{|\omega|}{2} (ab + ba) a
\]
calculating and extensively using the commutator relationship $ba = ab + u$ 
to move all of the ``b's'' to the right gives

\begin{gather*}
   [a,H] = \frac{|\omega|}{2}([a,ab] + [a,ba] ) \\
   [a,ab] = a^{2}b - a(ba)  \\
    [a,ab] = a^{2}b - a(ab+u)     \\
    [a,ab] = -au   \qquad  similarly \\
     [a,ba] = -au    \qquad hence \\
     [a,H] = -au|\omega |.
\end{gather*}
Thus
	\[ i\dot{a} = -\omega a \qquad therefore
\]

	\[ a_{s} = e^{i\omega s} a_{0}
\]
and we have the evolution equation for the operator ``a'' with respect to the time ``s'' which is of course the expected result from
standard quantum theory. Actually, we were expecting an additional $ \frac{\omega}{2} $ coming from the vacuum energy. We will be looking into this later.
A similar result for the operator "`b"' is

	\[ b_{s} = e^{-i\omega s} b_{0} .
\]

Clearly, with respect to the current set of macro observers, the above results show that there are two possible "`rest"' frames  for the operators "`a"' and "`b" and consequently "`H"'. Thus, if "`time"' follows the tracks given by the condition

	\[ u_{\omega} u_{s} = \pm 1
\]
where we use the notation $ s = |s|u_{s} $  with $u_{s} = e^{i \phi_{s}}$ as the corresponding phase factor for "`time", then we will have a "`rest"' condition. Hence we have defined a direction and orientation of "`time"'. 

If this condition does not hold, the macro observers would report exponential growth or decay and thus a possible unstable physical interpretation unless one is content with short time periods or unless the observers attempt to find a "`rest"' frame to perform calculations.

 Along these same lines and more formally, we prove that complex oscillators can be quantized and  that there is an inherent negative energy, anti-particle, time reversal and even the hint of spin 1/2  characteristics in quantum complex oscillators \cite{Car}, \cite{Fey3}.

\subsection{Complex Quantum Oscillator Theory; i.e., Weyl Algebra Theory}
As noted above, we have since learned that this section is essentially Weyl algebra theory. However, we want to keep the ideas and notations on record because it guides our thinking on the subject giving an almost totally symbolic/algebraic approach to quantum complex oscillator theory and possibly quantum field theory.

 Our presentation will appear to be almost the same as standard quantum mechanics and quantum field theory with the difference that we will be trying to stay as symbolic as possible and will try to be more careful in the notations, and of course, we want to directly include complex numbers into our calculations so that complex energy, complex time etc., make sense. 

 This is pretty much in the spirit of physicists' techniques. However, to perform actual calculations we will introduce a method to obtain quantitative results. Heisenberg and Dirac basically pioneered these techniques followed by Feynman and others. We need to work with non-unitary and non-hermitian operators to make progress in our work.

 To begin, we record a number of commutator calculations for the development of a symbolic and algebraic oscillator theory. 

\begin{definition}{}
Assume the following commutator condition
	\[ [b,a] = ba - ab = u
\] and define
	\[ H = \frac{1}{2}(ab + ba)
\]
where $ u $ is a phase times the unit operator. 
\end{definition}{}

\begin{definition}{}
Define the wave-particle algebra $ \Psi$ over the complex numbers with generators $"`a "'$ and $"'b"'$ denoted by
	\[ \Psi [a,b] 
\]with the commutator condition given above. This algebra contains the unit operator i.e., is a unital algebra.
\end{definition}{}

For the readers reference, this corresponds to the Weyl algebra $ A_{1}$ over the complex numbers.

 \begin{definition}{}
Define the left module 
	\[ \ominus = \Psi b
\] 
as the "`vacuum"'. 
We will be calculating modulo the vacuum and will denote this by $ \equiv $. Other vacuums, of course, can
be constructed such as those finite combinations of "`a's"' and "`b's"' or "`1"' that may be omitted in the module.
\end{definition}{}	
 \begin{definition}{}
Define the vacuum "`ket"' in the algebra $\Psi$ as follows

	\[ |0>^{\circ} = 1
\] 
then clearly we have the vacuum annihilation condition

	\[ b |0>^{\circ} \equiv 0 \qquad modulo \ominus 
\].
\end{definition}{}
\begin{definition}{}
(d)Define the following $"'a"'$ particle excited states lifted from the vacuum

	\[ a^{k} |0>^{\circ} = |k>^{\circ} = a^{k}
\] 
where in this context we will be favoring the production of "`a"' particles over "`b"' particles.
\end{definition}{}
\begin{definition}{}
(e)Define the conjugate of a "`ket"' which turns "`kets"' into "`bras"' where "`c"' is an arbitrary complex number:
	\[ \overline{ c|k>^{\circ}} = \overline{c}<k|^{\circ} =\overline{c} b^{k}.
\]
\end{definition}{}
\begin{definition}{}
Define the complex "`unitary"' operator to go from time "`0"' to time "`s"' by

	\[ U(s|0) = e^{-iHs}
\] 
so that on stationary kets with possible complex energy E modulo the vacuum $ \ominus = \Psi b$ at the time "`s=0"' with 
	\[ H |\psi,0>^{\circ} \equiv E |\psi,0>^{\circ}
\]
we have 
	\[ U(s|0)|\psi,0>^{\circ} \equiv e^{-iEs}|\psi,0>^{\circ} \equiv |\psi,s>^{\circ}.
\]
\end{definition}{}

Next in order to perform concrete computations we need to make contact with Schrodinger's wave mechanics and a semblance of Dirac delta function theory and distribution theory.
\par
Hence, when necessary identify 
	\[ a = u^{\frac{1}{2}}\frac{1}{\sqrt{2}}(p + iz) \stackrel{}{\rightarrow} u^{\frac{1}{2}}\frac{1}{\sqrt{2}}( -i\partial_{z} + iz)
\] with a similar identification for the "`b"'.

Note:  Later, we will see, that we can set the phase $ u = 1 $. Clearly, we have an ambiguous scaling of the commutator in that the phase factor u can slide into the factors different ways in order to normalize it to $ [b',a'] =1 $.  

We could have the balanced situation
\begin{gather*}
a' = u^{\frac{1}{2}} a \\
b' = u^{\frac{1}{2}} b  \\
\qquad or \\
a' = a \\
b' = \frac{b}{u}  
\end{gather*} and so forth.

 Our objective is to define the amplitude $ <z|k>^{\circ} $ for the variable  z.  However, this does not make sense, because z is not a label defined in the wave particle algebra $ \Psi[a,b ]$ and the previous definitions.
We get around this by using the following definition. 
\begin{definition}{}
Let $ z_{0}$ be a given fixed complex number.

Define the following matrix element 

	\[ < \widehat{z_{0}} |0>^{\circ} \equiv e^{-\frac{1}{2}z_{0}^{2}}
\]
extend this definition for arbitrary elements $\psi $ in the algebra $ \Psi[a,b]$ 

	\[ \widehat{<z_{0}}| \psi |0>^{\circ} \equiv \psi \widehat{<z_{0}} |0>^{\circ}
\]using 

	\[ \widehat{<z_{0}}| a^{k} |0>^{\circ} \equiv ( \frac{u^{\frac{1}{2}}}{\sqrt{2}}( -i\partial_{z} + iz))^{k} \widehat{<z_{0}}|0>^{\circ}.          
\]
\end{definition}{}
Then the following can be proved quite easily using standard quantum techniques.
\begin{lemma}
Using the definitions above:

(a) Commutator relations for "`a"' and "`b"',

\begin{gather*}
[b,a^{n}] = b a^{n} - a^{n}b = una^{n-1} \\
 [a,b^{n}] = a b^{n} -b^{n} a = -unb^{n-1}.\\
\end{gather*}

(b) Commutator for H

	\[ [H,a^{n}] =  una^{n} .
\]

(c) Norm of particle states modulo the vacuum $ \ominus = \Psi b$
	\[  <n| ^{\circ}n>^{\circ} \equiv b^{n}a^{n}  \equiv u^{n}n!
\]

(d) Eigen-energy of particle states modulo the vacuum $ \ominus = \Psi b$
	\[ H |k>^{\circ} \equiv u(k+\frac{1}{2}) k>^{\circ}
\]

\end{lemma}
Proof of the Lemma (highlights):
\par
Here we give some of the details for (c).
We can simplify the calculations in the following way:
absorb the "`u"' into say "`a"' and then define the new "`coordinates"' for the operators
\begin{gather*}
   a' = \frac{a}{u} \\
   b'=b \mbox{ thus} \\
   [b',a'] = 1 .
\end{gather*}
which in C* notation is $ [a,a^{*}]=1$,  hence all of these results should be quite familiar.
\par
Thus we can assume $ [b,a] =1$  and start the proof of (c).

By induction, assume the following for $ n $, where the case $n=1$ is obviously true modulo the vacuum $ \ominus = \Psi b$,

	\[b^{n}a^{n} \equiv n!
\] 
then using the results from (a) we have

\begin{gather*}
   b^{n+1} a^{n+1} = (b^{n+1} a )a^{n} = (ab^{n+1} + (n+1)b^{n} )a^{n} \qquad thus\\ 
   b^{n+1} a^{n+1} \equiv  (0 + (n+1)b^{n})a^{n} \quad hence \\
   b^{n+1} a^{n+1} \equiv (n+1)n! =(n+1)!  .
\end{gather*}
Note that the term 
	\[ b^{n+1} a^{n} \equiv b( b^{n} a^{n}) \equiv b n! \equiv 0
\] and we are done.
This concludes the proof of the lemma.
\\
 Using these results we can prove the following result:
\begin{theorem}
 (a) Complex quantum oscillator calculations can be conducted without introducing a "`vacuum"' per se or without introducing a metric structure using
bra and ket vectors in an analytical sense. These concepts are introduced algebraically.
\\
\\	
(b) All calculations can be carried out symbolically and algebraically, and thus convergence issues and renormalization issues can be performed more rigorously. "`States"' are just extended products and sums of "`a's"' and "`b's"' over the complex numbers.
\\
\\
(c) A particular orientation of "`time"' and vacuum state  will select the type of "`real'"'  particle to be observed.

\end{theorem}

 Next we proceed to an analytical method to quantize complex oscillators. This amounts to a Wick rotation method which pays closer attention to the Riemann surface aspects of the analytic continuation.

\subsection{Analytical Proof of the Quantization of Complex Oscillators}
\begin{theorem}

(a) The complex oscillator

	\[ H = H_{(u)} = \frac{1}{2}(p^{2} + u^{2}z^{2})
\]
can be put into a "`rest"' frame 

	\[H' = H'_{(1)} = \frac{1}{2}(p'^{2} + z'^{2})
\]
through scaling transformations  if and only if one chooses a phase w such that
	\[w^{4} = u^{2}
\]
and applying the following transformations

	\[z'= wz
\]
	\[p' = w^{-1}p
\]
	\[ H' = w^{-2} H
\]
	\[s' = w^{2}s  
\]
which maintains the quantum condition,

	\[[p,z] = [p',z'] = -i . 
\]
Only the case with $ w^{2} = -u $ leads to negative energy and time reversal

	\[ H' = -u^{-1}H
\] and
	\[ s' = -u s 
\] 
and the reversal of $a$  and $b$ to $ a = wb'$ and $ b = wa'$
where
	\[a = \frac{1}{\sqrt{2}}( p + iuz) 
\] and
\[b' = \frac{1}{\sqrt{2}}( p' - iz') .
\] 
(b) In particular for 

	\[ u = 1
\] there are two and only two possible rest frames with

	\[ w = \pm1
\]
and
\[ w = \pm i .
\]
Only the case with $ w = \pm i $ leads to negative energy and time reversal
	\[ H' = -H
\] and
	\[ s' = -s 
\] 
and the reversal of $a$  and $b$ to $ a = wb'$ and $ b = wa'$.
 Note: Quantum oscillators are inherently ambiguous with respect to a "`rest"' frame. This is not surprising since we are extracting
square roots using complex numbers and thus there will always be a "`+"' or "`-"' solution i.e., we will always have to introduce Riemann surface methods or explicit algebraic constructions to keep track of the ambiguities.

(c) Complex oscillators can be quantized using analytic continuation methods and staying in compact regions of time and space or essentially compact regions of time and space.
\end{theorem}
Proof:
(a): Using scaling transformations, the most general oscillator Hamiltonian with all parameters complex and all variables complex
	\[H = \frac{1}{2m}(p^{2} + m^{2}\omega^{2}z^{2})
\] and quantum condition
	\[    [p,z] = -i\hbar
\]
can be reduced to the form

	\[ H = \frac{1}{2}(p^{2} + u^{2}z^{2})
\] where "`u"' is a phase and 

	\[ [p,z] = -i .
\]
Choose a phase "`w"' such that
	\[ w^{4} = u^{2}
\] then scale as follows
	\[ z' = w z
\]
	\[ p' = w^{-1} p
\] then this transforms H and maintains the commutator condition
	\[ [p,z] = -i = [p',z']
\] as follows
\begin{gather*}
 H = \frac{1}{2}(p^{2} + w^{2}w^{2}z^{2}) \\
 H = \frac{1}{2}(p^{2}\frac{w^{2}}{w^{2}} + w^{2}w^{2} z^{2}) \\
  H =  \frac{1}{2}(p'^{2} w^{2} + w^{2} z'^{2} ) \\
  H = w^{2} H' \\
  H' = w^{-2} H .
\end{gather*}
To prove the swapping of the "`a"' and "`b"' for $ w^{2} = -u $ consider

\[a = \frac{1}{\sqrt{2}}( p + iuz) 
\] 
which leads to

\[a = \frac{1}{\sqrt{2}}( p + i(-w^{2})z) 
\] thus

\[a = \frac{1}{\sqrt{2}}( p\frac{w}{w} -iw (wz)) 
\] so
	\[ a = wb' .
\]
Note that this portion of the demonstration works in the opposite direction which essentially completes the "`iff"' proof.
Now time "`s"' transforms as follows using the Heisenberg evolution equations where "`v"' is a arbitrary dynamical variable from the algebra $\Psi[a,b]  $ over the complex numbers of the "`a"' and "`b"' operators
\begin{gather*}
    i \frac{d}{ds}v = [v,H]    \\
    i \frac{d}{ds}v = [v,w^{2}H'] \\
    i \frac{d}{dw^{2}s}v = [v,H']  \\
     s' = w^{2}s 
\end{gather*}
 and this finishes part (a). Clearly (b) is a special case of (a).
\\
Part(c) is straightforward.

\section{QFT and NCG Particle/Wave Analogies}
 Here we incorporate some of the motivation behind our investigation into RH using QFT ideas. We provide a ``particle-wave'' interpretation of Hasse's result where function fields correspond to "`waves"' with exponentials of "`s"' while the "`particles"' are created once they "`hit"' the vacuum from the corresponding operator. Wave particle duality is manifest because of the close association of the exponentials in $ s $ times creation operators in the corresponding field expressions.
 
 The idea of our creation operator paradigm came about in previous work on quantum "`anti-oscillators"' as we previously mentioned. We now know that Connes' NCG (Non-Commutative Geometry) \cite{C1}, \cite{JM}, \cite{CCM1}, \cite{CCM2} and others
is the way to proceed.
\subsection{A Quantum Eta Function Approach to RH}
In fact, we carry out some of the steps in Section 7. Hence as the writing of the paper proceeded, we have been gradually blending these ideas into our thinking and have incorporated some of the NCG language into our work i.e., cyclic homology concepts \cite{Lo}.
 
 We now proceed to develop the analogies. First we extend macro time-space to complex time-space
	\[ s^{\mu} = (s,z).
\]
Next we examine the following long loop integral along the branch cut considered in Theorem 2
	\[ \int_{\gamma(0)}\frac{z^{s-1}}{H_{n}(z)}dz \qquad ( 0 < s < n+1)   
\] 
which we can  call the $ H_{n}$, instead of a "`harmonic polynomial'', a "`cyclic  polynomial
	\[	H_{n}(z) = (z+1)(z+2)...(z+n+1).
\]
It should be clear from the discussion in Sections 1-3, that this integral and associated integral have a direct bearing on RH. 

 Next we apply our original concept of quantizing complex variable function theory using creation operators i.e., QFT to "`soften up"' the analysis so that we can study more delicate "`phase"' locking \cite{PH} issues in the alignment of the non-trivial zeros of the $ \zeta(s)$. The "`soften up"' idea is one of the central notions of NCG.

Thus using Weyl creation operator correspondence that we have mentioned several times previously
	\[ z + ik \stackrel{}{\rightarrow} p + ikz/\sqrt{2} = a_{(k)} 
\] 
 we can perform symbolic particle quantum field theory.
 
 The reader can find an excellent introduction to these concepts in Dirac \cite{D2} and many other books e.g. \cite{K}.  
We merely mention that there are bosonic commutator relationships and fermion commutator relationships such that if the operators commute, they are of the same type or they belong to different states. 

And if there are commutators involving creation and destruction operators in the same state, we obtain a non-zero result with Dirac deltas times $\hbar$. So what we are doing is "`second quantizing"' complex variable theory. We have also considered setting up a Feynman path integral formulation as well \cite{Fey2}. However, we made more progress with the particle methods introduced here.

 The next idea is to replace the above integral into "`new"' coordinates using these operators to obtain something that looks like the following

\[ \int\frac{a_{0}^{s}}{a_{0} a_{1} a_{2}...a_{n}a_{n+1}}da_{0}   .
\] 
Hence, we have a highly interacting finite particle QFT type of analysis to consider. And  we have the task to put all of this on a rigorous foundation. Our paper started some of this; however, NCG seems to have already created most of the tools for moving forward in this direction. So our strategy is to combine NCG methods with our own methods to achieve the desired results. 

Now Theorem 2 shows there should be a relationship of the above integral to a finite approximate Hasse eta function i.e., local eta

	\[ \zeta^{H}_{n}(s) = \eta_{n}(s) =\sum^{n}_{k=0} (-1)^{(k)} \binom{n}{k}(k+1)^{-s} 
\]
which we will now interpret using  concepts as foundation from Theorem 3 and Theorem 4. Recall that Theorem 3 gives us the right to pursue symbolic/algebraic methods in our computations. However, we only sketch the details in this paper. 

The above operator integral is replaced by the more "`gentle"' quantum integral

	\[ \eta_{n}^{q} (s)   = \int_{\gamma} a_{0}^{s} b_{0} b_{1} ... b_{n}b_{n+1}da_{0} .
\]
We then introduce some cyclic homology concepts from Loday \cite{Lo}
and call the integrand the target complex 

	\[ K_{n} = a_{0}^{s} b_{0} b_{1} ... b_{n}b_{n+1}
\]
which is to be analyzed using insertion and deletion operators with respect to a well defined algebra of appropriately defined a's and b's. 
We then define the sampling complex as follows

	\[ \gamma_{n} = [a_{0},a_{1},...,a_{n},a_{n+1} ]
\] 
with the following preliminary definition of the differential $ d a_{0} $ at a sample point $ a_{k} $ as

	\[ da_{0} = a_{k}.
\]
Then we proceed to define integration of these operators, with respect to an appropriately defined vacuum, as follows
	\[ \gamma_{n} K_{n} =  (a_{0}da_{0}) K_{n} + (a_{1}da_{0}) K_{n} +..+ (a_{n}da_{0}) K_{n} + (a_{n+1}da_{0}) K_{n} .
\]
The time-space thread or "`micro-quantum arrow of time"' is chosen to be as tight as possible energetically speaking, allowing only quantum fluctuations not to exceed, for
either the a's or b's, the absolute energy of 

	\[ k+1 +\frac{1}{2}   .
\]
at a given time step "$k=-n-1,-n,..-1,0,1,2,..,n,n+1$"'. Positive energy is equivalent to the number of particles in our discussion here \cite{Wyl2}.

The a's and b's are decomposed further into the form

	\[  a_{k} = \frac{a^{k+1}}{\sqrt{(k+1)!}}
\]
and
 \[  b_{k} = \frac{b^{k+1}}{\sqrt{(k+1)!}}
\] 
with the usual commutator, as before,

	\[ [b,a] = 1.
\]The given macro"arrow of time"' is the standard one given by the complex
	\[ [b_{n+1},b_{n},...,b_{0},a_{0},a_{1},..,a_{n},a_{n+1}]
\]which is to be used in carrying out micro-quantum calculations and what we call the "time-space thread"' with respect to a well defined vacuum very similar to the above complex but expressed in module form. Other micro-quantum time-space threads are possible.

For example the vacuum,
	\[ \ominus_{n} = \Psi b_{n+1} +\Psi b_{n}+...+\Psi b_{0} + a_{0}\Psi +...+ a_{n}\Psi+a_{n+1}\Psi
\]
has exactly $ 2^{n+1}$ typical "`events"' stored here in the form
	\[[ (a_{0} "`or"'  b_{0}),...,(a_{n} "`or"'  b_{n} )  ]  
\]
in the forward micro time steps $[0,1,..,n]$. We do not take $(a_{n+1} "`or"'  b_{n+1} ) $ as one of our typical events because it is at $ \infty$, c.f. Figure 2. The time-space thread runs vertically upward.
	\[
\]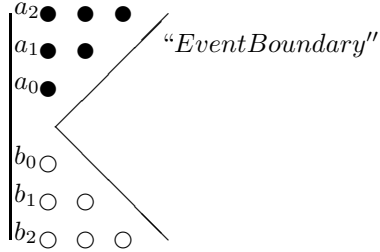
\begin{figure}[h]
\caption{ The time-space thread for $ \ominus_{1}$. }
\centering 
\setlength{\unitlength}{1cm}
\begin{picture}(2,2)
\put(.6,0){\line(1,1){1.5}}
\put(.6,0){\line(1,-1){1.5}}
\put(0,0){\line(0,-1){1.5}}
\put(0,0){\line(0,1){1.5}}
\put(.5,.5){\circle*{.2}}
\put(.5,1.0){\circle*{.2}}
\put(1,1.0){\circle*{.2}}
\put(.5,1.5){\circle*{.2}}
\put(1.0,1.5){\circle*{.2}}
\put(1.5,1.5){\circle*{.2}}
\put(.5,-.5){\circle{.2}}
\put(.5,-1.0){\circle{.2}}
\put(1,-1.0){\circle{.2}}
\put(.5,-1.5){\circle{.2}}
\put(1.0,-1.5){\circle{.2}}
\put(1.5,-1.5){\circle{.2}}
\put(.05,.5){$ a_{0} $}
\put(.05,1.0){$ a_{1} $}
\put(.05,1.5){$ a_{2} $}
\put(.05,-.5){$ b_{0} $}
\put(.05,-1.0){$ b_{1} $}
\put(.05,-1.5){$ b_{2} $}
\put(2.,1.){$ ``Event Boundary''$}
\end{picture}
\label{fig:time-space} 
\end{figure}
	\[
\]
	\[
\]
	\[
\]
The b's will interact with the a's in the upper thread.
Notice that we will "'induce"' the b's to sit on a given forward time step and that we do not allow
two different b's or two different a's to be there at that same 'time"' for very long. Then we process all of these events and watch how they drop from vacuum to vacuum recording useful byproduct information,, i.e., these are the symbolic "`photons", until we reach the scalars or a useful stopping point for further analysis. We plan on performing some sort of ``quantum particle fluid'' analysis.

 In the above algebra, we also use the peculiar quantum projective concept that $ "`zero"' $ and $ "`\infty "`$ are equivalent, in some sense, as follows
	\[ a_{n+1} \equiv b_{n+1} \equiv \infty \equiv 0 \equiv a_{0} \equiv b_{0}.
\]This holds because of the design of our vacuum. 

Next, with these definitions, we plan on studying the following matrix element, which, we believe, is analogous to cyclic cohomology constructs,
	\[ <0,a_{0}^{s}|\eta_{n}^{q} (s)> .
\]
Note that a variety of additional correlations can also be studied as well.
The above will then be related to Hasse's local eta of the form

	\[ \eta_{n} (-s) = 1 -\binom{n}{2}2^{s} + \binom{n}{3}3^{s} -....  .
\]

 Thus we interpret the local etas as relative complex amplitudes and the global eta as the superposition of the complex amplitudes of the local etas weighted by the number of quantum states in the corresponding finite QFT.

The next step is to study micro-equilibrium conditions near $ s = \frac{1}{2} + it$,with t real, and prove that the micro-equilibrium is achieved for well defined t's with only very small decay $id$ in t as follows:

	\[ t = t_r + id
\]
where $t_{r}$ is real and $d$ is real i.e. these are the "`proto-zeros"', we alluded to before.

 In general the proto-zeros, we expect can be found, for a fixed $ \sigma $, for variable "`t"' in $ s= \sigma + it $ using simultaneous Diophantine approximation methods \cite{HW} over the rationals by looking at the ratios of the logarithms of the primes, which will be transcendentals \cite{S}, \cite{Did3}, \cite{V} in the corresponding interval for a given local eta. There will be an infinite number of "`good"' "`solutions". Hence there will be an infinite number of proto-zeros at any given local eta. The corresponding proto-zeros for large $n$ for a given $s = \sigma +it$ should be good approximations to a global zero as well \cite{B1}. 

 In carrying out these computations we found it convenient to work with a local Planck's constant \cite{L2} defined by

	\[\hbar_{p} = \frac{log(p)}{2\pi}.
\]
which will measure the "`resolving"' power of a given local eta's ability to detect the largest
prime in its corresponding lattice. The smallest oscillatory time interval a local eta can resolve is given by

	\[ \frac{1}{\hbar_{p}}
\]
 Hence a proto-zero for a given local eta, at a given time-space point, is when Hasses' local eta is as small as possible relative to a small domain containing the point.

 The meaning that we attach to a typical proto-zero is that the local quantum prime counting sieving operator is delivering its assessment in the form of a complex amplitude near the space time point $s = \sigma +it $ that it is "`satisfied"', within its resolving power, in counting the number of primes less than or equal to $t$. A non-zero complex amplitude, with "'decay", carries important information useful to the "`collective"' in determining a global zero. So "`iso-phase"'  \cite{Der} curves may be a more general concept to consider. This now leads us to discussing some general technical strategies for settling RH. 

\section{ Overview of Technical Approaches to Settle RH}

\textbf{Method 1: Classical with Operator/QFT and NCG Inspired Insights (Local to Global and Global to Local)}

The idea is to prove that the "`cloud"' or "`collective"' of local proto-zeros
helps to determine the global zeros of $\eta(s)$ and conversely that the global zero is the "`center of gravity"' of the local zeros.

 So if the cloud of local zeros moves then the corresponding global zero must move maintaining its center of gravity. We intend to prove that any given global zero is determined by only a finite number of proto-zeros. We think this can be handled by a suitable fixed point theorem.

 We expect that the latter set of proto-zeros will not be exactly aligned on the critical line and we expect some small "`decay"' in each one will be necessary to achieve this representation. We expect this because a given local proto-zero by itself is "`too
loud"' and needs to be "`detuned"' for the other local etas to have their say. 

 We also intend to prove that the local etas are in micro-equilibrium on or very near the critical line. So, if a global zero would not be on the critical line, it would move the cloud into a disequilibrium state contradicting our previous "`result"'. Note that Conrey reports that 99 per cent
of the global zeros are very close to the critical line \cite{Con2}. In fact, Conrey proved that at least 40 per cent of the global zeros are on the critical line \cite{Con1}.

 Clearly with these steps rigorously established, we will have settled RH.

 Possible exceptions would be handled by iterative bootstrapping methods familiar to functional equation theorists \cite{A}.

 We will use the operator NCG and QFT insights above coupled with the classical methods working directly with the local eta functions. Hence, we expect to obtain rigorous results independently of NCG or QFT.  

 Note:  Our Theorems 1 and 2 gave some important clues into the structure of the local etas, but the reader is undoubtedly aware that Hasse's result:
	\[ \eta(s) = \frac{\eta_{0}(s)}{2} + \frac{\eta_{1}(s)}{2^{2}} + \frac{\eta_{2}(s)}{2^{3}} + ....
\]can be proved using elementary methods with the Binomial Transform and that this amounts to a renormalization of $\eta(s).$
\\
Therefore, we plan on staying, in this approach, as elementary as possible in the proofs.
\\

On the other hand, in a complementary fashion, using similar concepts:
\\
\textbf{Method 1': NCG Methods (Local to Global and Global to Local)}
We would try to obtain rigorous results using the NCG/QFT concepts that we outlined in this paper while implementing Method 1 above in parallel.  We would try to utilize the most advanced methods that we can understand using Connes' and his collaborators' techniques or others' methods e.g. Conrey's 
and his collaborator's methods.
\\
\\
\textbf{Method(s) 2: NCG Methods Applied to Previous Techniques}
For example, we would try to obtain local theta functions relationships using the ideas of this paper and re-examine Hardy's method etc.,
including the author's slow continued fractions method in the fundamental equation of information theory using ergodic theory \cite{Did1}.

 The following section now outlines such an approach with ergodic theory and Euler factors.
 
\section{Approach to RH via Euler Factors}

We sketch a Euler Factor approach to RH using Weyl algebra creation operators to analyze the quantum behavior of the integers as one pushes up the critical line near a macro zeta function zero. This research stems from the eta functions\cite{Did2} methods described in Sections 1-3 with Ornstein's \cite{O},\cite{O1} gadget constructions from ergodic theory \cite{Wa} plus the author's \cite{Did1} previous research concerning the fundamental equation of information theory\cite{AD}. Maxim and his collaborators \cite{M1} results provided impetus to move forward in this direction.

 We intend to show that there is a quantum $ \chi_{q} $ symmetry, where notationally, q is a label denoting "`quantum", such that
	\[ b_{(1)}^{s} a_{(1)} \zeta_{q}(s) \equiv \zeta_{q}(s)
\] which causes the zero to snap onto $\sigma = \frac{1}{2}$ because
	\[\bar{\zeta}_{q}(s) \zeta_{q}(s) \equiv 1 + 2s(s-1) + 2s(s-1) \equiv \eta(s) .
\] and by algebraic coherency \cite{S1}; by ergodic flow construction; and the presence of three renormalized quantum fields: "`0"',"'1"',"'-1"' with appropriate phases \cite{Bin}.

 We expect deeper results can be obtained by this process coupled and contrasted with known methods e.g. Nakajima creation operator techniques for Douady spaces \cite{CM} and Hilbert schemes \cite{M1} and the recent results of Ngo \cite{N}. 
 
\subsection{Introduction of Weyl Algebra Creation Operators: One for Each Non-Negative Integer}
 
 First recall the classical Euler's result
	\[ \zeta(s) = 1+ \frac{1}{2^{s}} + \frac{1}{3^{s}} + \frac{1}{4^{s}} +.. =  \prod_{ p = primes }\frac{1}{1-\frac{1}{p^{}s}    } \qquad ( \sigma > 1).
\] where $ s= \sigma +it$ \cite{E}. Next we need a quick review of Weyl algebras.

As noted in Sections 1-3,we introduce the corresponding Weyl algebras of commutators\cite{Hart} 
\begin{gather*}
 [b_{(k)},a_{(k)}] = 1 \\
 [b_{(k)},a_{(k')}] = 0 \qquad (k \neq k') \\
 [b_{(k)},b_{(k')}] = 0 \\
 [a_{(k)},a_{(k')}] = 0 \\
\end{gather*}
where the Weyl algebra is built up algebraically from a free algebra on the symbols corresponding to the a's and b's  and adjoining the two sided ideal generated by the above relations. 

 For example, using our previous notations \cite{Did1}, the particle-wave algebra on a and b is given by free algebra over the complex numbers modulo the two sided ideal defined by the commutator condition

	\[ \Psi = \frac{\mathbb{C}[a,b]}{( [b,a] - 1)} = A_{1}(\mathbf{C}).
\]
We have waves in particle-wave algebra $\Psi$ because of the $``up''$ and $``down''$ characteristics of the $a's$ and $b's$ with the complex exponentials giving the macro time-space information, c.f. the next section.

\textbf{Note}: We rename these operators Weyl creation operators because of their close association with Quantum Harmonic Oscillator theory. 

\subsubsection{Review of Complex Quantum Oscillator Theory}
 It is useful to repeat some of the previous analysis on complex quantum oscillator theory in this context.
 
 The creation/destruction operators are defined as follows, using our paradigm of replacing a linear polynomial over the complex numbers into a creation operator as follows
 
	\[ z + ik \stackrel{}{\rightarrow} \frac{1}{\sqrt{2}}( p + ikz) = a_{(k)}= a
\]  and with
\[ z  \stackrel{}{\rightarrow}  a_{(0)}= 1
\]
where on the left side $z$ is a complex variable; k corresponds to complex number i.e. a scalar with phase u and $k = |k| u$; and on the right side $p$ now corresponds to a momentum operator in the Weyl algebra and is not to be confused with primes $p$; and $z$ now corresponds to a position operator in the Weyl algebra.
  Also we have
\[  a_{(-k)}= b .
\]
 Using the standard commutator condition CCR (Canonical Commutator Relations) \cite{Petz}
	\[ [p,z] = pz-zp = -i
\] and the corresponding quantum harmonic oscillator mathematics \cite{Bo} leads to the condition with these definitions that
the symmetrized product corresponds to an energy of $k = |k| u$  where $ u = e^{i\theta } $ because

\begin{gather*}
H_{(k)} = H= \frac{1}{2}(ab +ba) \\
ab = H - \frac{1}{2}k  \mbox{ the energy goes down by k/2}\\
ba = H + \frac{1}{2}k  \mbox{ the energy goes up by k/2}  \\
H = \frac{1}{2}(ab +ba) \\
[b,a] = ba-ab = k
\end{gather*} as one can verify by multiplying the corresponding products in our definitions above. This is usually expressed in physics as 

	\[ [a,a^{*}] = 1
\] where $a=b $ and $a^{*} = a$ in our notation.

 This shows that we will be scaling and normalizing all of our creation operators definitions to satisfy
	\[ [b,a] = 1.
\] For example
\begin{gather*}
 a \stackrel{}{\rightarrow} u^{{-\frac{1}{2}}}|k|^{-\frac{1}{2}} a \\
 b \stackrel{}{\rightarrow} u^{{-\frac{1}{2}}}|k|^{-\frac{1}{2}} b.
\end{gather*}
 If we want to find out the detailed energy, momentum, and position of a given state, we will need to take that into account.
 \par
 \textbf{Remark :} \textit{Thus an important intuitive concept here is that we create intrinsic particles and that we attach energy and other particles in our information-theoretic constructions to it by which we transport information locally and globally as needed. Further the b's correspond to gentle multiplicative inverses i.e. $ ba \equiv 1$ .}
\\
 So the reader should be aware that all of our creation operators have been scaled and normalized and that the ambiguity in the square roots have been accounted for. 
 \par
\textbf{Note}: What we have done is to write all the operators with the individual momentum and position joined together and call these operators creation operators. 
 \subsection{Replace the non-negative integers with creation operators}
Now in the classical Euler factor expression given above, replace with Weyl creation operators, one for each  non-negative integer $m$, to obtain a quantum zeta and quantum Euler factors in the critical strip 

	\[ \zeta_{q}(s) = a_{(0)}+ a^{s}_{(1)} +  a^{s}_{(2)}+a^{s}_{(3)} +a^{s}_{(4)}.. \equiv   \prod_{ p = primes } E_{(p)}^{q}(s)
\]
where
\begin{gather*}
 a_{(0)}=1  \\
 E_{(p)}^{q}(s) = a_{(0)}+a^{s}_{(1)} +  a^{s}_{(p)}+a^{2s}_{(p)} + a^{3s}_{(p)}.... .
\end{gather*}

This holds almost formally in an abstract ring of power series sense. It gives the arithmetic blending instructions needed later in our constructions.
Mathematical convergence and other issues will be addressed later. The analysis is quite delicate and subtle. The notation $\equiv$ will be explained later. It is an algebraic geometry condition analogous to the calculation of the cotangent vectors in algebraic geometry i.e. $M/M^{2}$ \cite{V},\cite{Miln1}.
\\
 Hence, we are automatically in a Noncommutative Geometry\cite{L} setting with the corresponding complications in the analysis. 
 However, we now are in a very advantageous situation because we have created particles and waves from which we can watch and mold very delicate arithmetic, topological and algebraic operations all together. The waves allow us to extend beyond classical boundaries and to have the ability to be in many places at the same time. 
\par 
  So the basic idea is to quantum mechanically transport the real line together with the integers onto the critical line using quantum field theory concepts.

 In order for us to accomplish this, we will need to calculate a complex power of a creation operator $ a^{s}$ without using logarithms.
The next section briefly outlines the methodology for this calculation.

\section{ Calculation of $a^{s}$ }
 Page 1 of Janos Kollar's book \cite{K} on resolving singularities gives the clue how to proceed. The binomial expansion leads to 
	\[ a^{s} =(1 + (a-1))^{s}= \sum^{\infty}_{k=0} (a-1)^{k} \binom {s}{k}
\]
 where
	\[\binom {s}{k} = \frac{(s)}{k} \frac{(s-1)}{1}.....\frac{(s-(k-1))}{k-1}.
\]
Also observe that 
	\[\frac{s-m}{m} = (-1)(1-\frac{s}{m})
\]
hence
\[\binom {s}{k} = (-1)^{k-1}(\frac{s}{k}) ( 1-\frac{s}{1}).....(1-\frac{s}{k-1}).
\]
\\
Thus we have injected, through inclusion, the algebra based on the $a$ and $b$ into a power series in $a'=a-1$ and $b'=b-1$. Observe that
	\[ [b-1,a-1] = [b',a'] = 1
\] 
which we think of as translation type or transitional creation operators \cite{BB}. This operator defines a direction of "`time"' used by the observers. It is not difficult to demonstrate that the energy of these operators, using symmetrized products modulo the vacuum $\ominus_{0} = \Psi b + a\Psi $ is one more than the energy of the original a and b confirming our intuition. One can easily see this also because $ -1 = i i $ providing that $z$ is near the vacuum $1$. Thus we have renormalized $a^{s}$.
\\
\begin{definition}{}
Next define the following "`quantum coherence"'  $ | |_{q}  $functional as follows:

	\[ |\sum_{k=0}^{\infty} a'^{k}\binom {s}{k}|_{q} = |\sum_{k=0}^{\infty}\binom {s}{k}|
\]
where we have used the discrete topology functional on the operators and the factors involving the scalars are analyzed using the classical norm defined on the complex numbers. 

\end{definition}{}
 This idea came from Serre's treatise on coherent algebraic sheaf theory \cite{S1}. One can systematically extend this notion to more general arithmetic, topological and algebraic situations. 
 Next introduce the function
	\[ \pi(s) = \sum_{k=0}^{\infty}\binom {s}{k}
\]
not to be confused with the $\pi(x)$ prime counting function. The convergence behavior of this function will determine how to extend $a^{s}$ locally using a theorem of the following type.
\begin{theorem}
The complex power of $a^{s}$ converges in the above sense if
\\
\begin{align*}
( \sigma-1)^{2} + t^{2} <  1 \\
\frac{1}{4} \leq \sigma \leq \frac{3}{4} \\
0 \leq t \leq \frac{1}{2}	\\
\end{align*}
 We use the notation $ \textit{F}_{a} $ to denote the above domain and the companion domain $ \textit{F}_{b} $. We call these domains Clifford Domains

\end{theorem}

A proof (sketch) utilizes the fact that the complex norm, for positive integers $k$,
	\[ |1-\frac{s}{m}|^{1/2}  < 1 \qquad { m=1,2,..,k-1}.
\]
and that any accumulating phases cannot produce a divergent harmonic series; otherwise we can perform a rotation back onto the real line with a rotated Clifford domain with its own critical strip which overlaps the original domain producing a contradiction since it holds in the real case. Only finitely many corrections are needed to handle possible divergences.
To start the proof, begin with finite 2-adic representations of $s = \sigma +it $ for example of the form
\[ s = \frac{1}{2} -\frac{1}{2^{3}} +i( \frac{1}{2} - \frac{1}{2^{5}} )
\]
to prove that there will be only a finite number of alternating type harmonic series(etas)with given "`locked"' phases i.e. a finite number of convergent etas. By using a large enough power of 2, we can make all of the phases to lie in a two dimensional lattice over the complex integers.  We are exploiting the fact that the result holds for the perfect phase cases e.g. $ s= i\sigma$ and with the angles 
	\[ \theta = \frac{2\pi}{2^{2n}}
\] finish the proof using the local polynomial topology on the overlapped domains.

To see a connection with Weyl algebra and Hilbert scheme methods consider, for a given $s$ at a given time step $k$, the operators
\begin{gather*}
a_{(m)} =  1- \frac{s}{m} \qquad  ( m = 1,2..,k-1)  \\
b_{(k)} = \frac{s}{m}  \qquad  thus\\ 
\binom {s}{k} = (-1)^{k-1} b_{(k)}a_{(1)}....a_{(k-1)} \qquad with \\
a_{(m)}a_{(n)} = a_{(n)} + a_{(m)} + (-1)a_{(mn)} \qquad furthermore \\
a_{(p)}a_{(p)} = a^{2}_{(p)}  \equiv (-1)a_{(p^{2})}
\end{gather*}
Further study of this situation helps to understand the Weyl creation operator method that we have been developing.
\par
To calculate $a^{i}$ use
\begin{align*}
a^{i/2} =  b^{1/4}a^{1/4+i/2}   \\
a^{i}	= a^{i/2}a^{i/2}.
\end{align*}
\\ These calculations will help to decide when there are suitable stopping points in the intermediate steps in the next constructions below.

\section{ Fusing (Folding) the Clifford Domains Up the Critical Strip: Squeeze, Stretch and Tile then Blend.}
First observe that the domains are one half of quantum in the real direction and one half quantum in the imaginary direction. Hence, have the same information content i.e. an area of one quarter. So the plan is to maintain a constant area of one quarter in all of our constructions. This is a basic concept in ergodic theory \cite{Sh} and information theory \cite{G}.

Now the idea is to construct the quantum symmetry $\chi_{q}$ which is analogous to Riemann's reflection functional equation  symmetry $ s -> 1-s $ defined as follows
\begin{align*}
\chi_{q} \textit{F}_{a}  = b^{s}a\textit{F}_{a}    \\
\chi_{q} \textit{F}_{b}  = a^{s}b\textit{F}_{b}
\end{align*}
 This is the fold operation. Clearly this operation is fusing the two domains together. 
\begin{itemize}
	\item{Step 1 Fuse and fold the domains together.} 
	\item{Step 2 Stop when the quantum symmetry stabilizes. }
	\item{Step 3 Squeeze the domains by 1/2 in the $\sigma$-direction toward the critical line. For example,
	the mapping, interpreted into a and b language, $ \sigma -> \frac{3}{2} \sigma -\frac{1}{4}  $ squeezes the left half domain by 1/2. The right domain is carried out in a similar fashion using the reflection symmetry.}
  \item{Step 4 Stretch the domains in the t-direction by similar transformations. Note: The above preserves the information content of a given phase cell.} 
  \item{Step 5 Tile to fill in the rest to obtain the expanded domains.} 
  \item{Step 6 Repeat as necessary to achieve the desired degree of accuracy.}
  \item{Step 6.1 These operations effectively extend the operators $a^{s}$ and $b^{s}$ up the critical line and 
  effectively transport the quantum symmetry $\chi_{q}$ along with it.}
  \item{Step 6.2 Furthermore, we require a s-fractional permutational scaling such that
\begin{align*}
b^{Ns} \stackrel{}{\rightarrow} b^{Ns} (N+1!)^{-s/2} = b_{(N)}^{s}	\qquad{and}\\
b^{s} a^{1} a_{(m)}^{s} = m^{-s} a_{(m)}^{s}
\end{align*}
  
  }  
  \item{Step 7 With this accomplished we study the behavior of $\zeta_{q}(s)$ near a macro zero to show that
	\[ b_{(1)}^{s} a_{(1)} \zeta_{q}(s) \equiv \zeta_{q}(s). \] 
 These operations are carried out for all the operators, i.e., for each integer and the unit translation operator $ a' = a_{(1)} - 1$,which gives the inherent quantum direction of "`time"', with additional blending carried out to properly express arithmetic information between the primes and the integers.   } 
	                                 	
\end{itemize}  
 We believe the above operations correspond to some sort of fractional cohomological extension theory similar to the theory of Hilbert schemes \cite{M1}. Clearly, we are using ergodic theory concepts here to transfer information to the particles.
 
\section{ Calculation of the products $a_{(p^{e_{1}}_{1})}^{s} ...a_{(p^{e_{n}}_{n})}^{s} $ }
Having made the above preparations, we expect to prove that the deviation for an arithmetic product near a macro zero is
	\[  d_{m} = a_{(p^{e_{1}}_{1})}^{s} ...a_{(p^{e_{n}}_{n})}^{s} -a^{s}_{(m)} \equiv 0
\]where
	\[ m = p_{1}^{e_{1}}..p_{n}^{e_{n}}
\] and where
the above set of primes correspond to the set of participating primes accounted for by the standard Riemann zeta theory less than or equal to $t$
for $s = \sigma + it $ which we think corresponds to s-fractional strata from Hilbert scheme theory.
The equivalence $\equiv$ corresponds to the observer modules of the form $ M=\Psi b + a \Psi $. 

\subsection{ Equilibrium Calculation}

Next we expect to  to find the "`trace"' of $ \zeta_{q}(s)$ as follows
\[\bar{\zeta}_{q}(s) \zeta_{q}(s) \equiv 1 + 2s(s-1) +2s(s-1)\equiv \eta(s)
\]
where the notation $ \bar{\zeta}_{q} $ means writing the expression in reverse order changing a's into b's and taking ordinary multiplicative inverses of scalars since, in the "`renormalized"'," thickened"' or "`dressed"' form , $\zeta_{q}(s)$ decomposes into three quantum fields of appropriate phases

	\[ \zeta_{q}(s) \equiv  \tilde{a}^{s}_{(0)}+ \tilde{a}^{s}_{(1)} +\tilde{b}^{s}_{(-1)}.
\]

So what ever symmetry was found in the former quantum field $\tilde{a}^{s}_{(1)}$ must occur in the latter quantum field $\tilde{b}^{s}_{(-1)}$ with reverse symmetry to maintain zero angular momentum. The pole field $\tilde{a}^{s}_{(0)}$ has been smoothed. Any back reactions will have been coherently stopped in our constructions.

This, if true, explains why the zeros are on the critical line.
\\
\\
 To see the possible truth of this, observe the following calculation with the cubic and higher terms ignored and then analyzed by the modules of the form $ M = \ominus_{0} = \Psi b + a\Psi $.
\begin{align*}
b^{s}	a^{s} = (1+ s(b-1) )( 1+ s(b-1) )  \\
b^{s}	a^{s} = 1+ s(b-1)+ a(b-1) + s^{2}(b-1)(a-1) \\
b^{s}	a^{s} \equiv 1 -2s + s^{2}(ba -a -b +1)\\
b^{s}	a^{s} \equiv 1 -2s + s^{2}(ba  +1)\\
\end{align*}
 Hence we arrive at
\begin{align*}
b^{s}	a^{s} \equiv 1 -2s + s^{2}(ba  +1)\\
b^{s}	a^{s} \equiv 1 -2s + s^{2}(ab +1 +1)\\
b^{s}	a^{s} \equiv 1 -2s + 2s^{2}\\	
b^{s}	a^{s} \equiv 1 + 2s (s-1)
\end{align*}

Also the observer modules $ M= \ominus_{0} = \Psi b + a\Psi $ have to be adjusted to account for the order of processing of these steps i.e. "`chirality"' \cite{BD} since we have to get rid of the additional "`1's"' that are produced in two of the steps that are being forced to be synchronized near the macro zero.

\section{ Concluding Remarks}
We believe that we have provided a second approach to RH, using a combination of classical thinking with QFT\cite{Sch} insights, with specific targeted potential theorems that should explain why RH is true. The first approach was outlined in Section 5 using the etas. Actually both approaches are linked together.

 To further aid the reader's intuition, it is helpful to think of the operators $ a_{(k)}$ as infinitesimals with the quantum ability to be in two places at the same time i.e., between neighboring integers additively and multiplicatively. It is helpful to think of the operators being decomposable into the form 
	\[ a_{(k)} \approx  a_{(1)}^{k}.
\]
and that the translation operator
	\[ a' = a_{(1)} -1
\]
helps to extend the operators into a "`fractional"' state in the critical strip and beyond and provides the micro direction of time. So this is similar to adele ring theory concepts in number theory \cite{We},\cite{Miln2} except that the integers can "talk" to each other beyond quadratic reciprocity.

 Further, it is useful to think of $ a^{s}_{(k)}$ as a " quantum field"' associated with the integer $k$ at a given complex time $s$; and that the algebra being constructed is the quantum field being analyzed in a given time-space neighborhood i.e., we are renormalizing or thickening the operators.  It helps to think that we are working along the sum of filtered products of a given fixed time-space topological \cite{Pa} space corresponding to a product of two complex projective manifolds modulo at divisible group times a two dimensional complex lattice. It helps to think that certain fluctuations from a "`vacuum"' are being looked at very carefully by our observers i.e. lower order operators have precedence over quantum fluctuations. Also the observers must change progressively as the constructions proceed.
 
   The concept of "`state-measurement"' occurs only when the corresponding operator "`touches"' the corresponding particles; one does not touch a particle until it is advantageous to do so. Further, symmetrized products of creation operators are the Hamiltonians that generate the fields.
\\
\\   
Further sources of information include: Kahler Geometry \cite{Mo},\cite{Gol}; Elliptic Functions \cite{Ha},\cite{Ko}; Algebraic Topology \cite{Sp}; Functional Analysis \cite{Shil}; Information Theory \cite{G},\cite{Go}; Dynamical Systems \cite{Mi},\cite{LF}; Cohomology \cite{St}; L-functions and Langlands Program \cite{Co}; Cryptanalysis \cite{Kopp}. 
 
 \textbf{Note} Our citations here do not necessarily imply competency by the author. This will be established through our proofs as the research is conducted..
   
\section{Acknowledgments}
Overall, I would like to acknowledge Eli Lubkin for the many conversations dealing with physics over the years which has inspired some of my research activities. And I would like to acknowledge some very helpful recent conversations with L. Maxim and J.Cogdell dealing with the Euler factor approach to RH.

\end{document}